\documentclass[12pt]{article}
\usepackage{amsmath,amssymb,amsthm,amsfonts,graphicx,latexsym,pdfsync}
\usepackage{mathrsfs,color}
\usepackage[latin1]{inputenc}
\numberwithin{equation}{section}

\DeclareMathOperator{\dist}{dist}

\DeclareMathOperator{\spt}{spt}

\newcommand{\R}{\mathbb R}

\newcommand{\be}{\begin{equation}}
\newcommand{\ee}{\end{equation}}
\newcommand{\ds}{\displaystyle}
\newcommand\eps{\varepsilon}

\theoremstyle{plain}
\newtheorem{thm}{Theorem}[section]

\newtheorem{lemm}[thm]{Lemma}

\theoremstyle{remark}

\newtheorem{exam}[thm]{Example}

\newtheorem*{ack}{Acknowledgements}
\newcommand{\bib}[4]{\bibitem{#1}{\sc#2: }{\it#3. }{#4.}}

\title{Optimal location problems with routing cost}

\author{Giuseppe Buttazzo, Serena Guarino Lo Bianco, Fabrizio Oliviero}

\begin{document}

\maketitle


\begin{abstract}
In the paper a model problem for the location of a given number $N$ of points in a given region $\Omega$ and with a given resources density $\rho(x)$ is considered. The main difference between the usual location problems and the present one is that in addition to the location cost an extra {\it routing cost} is considered, that takes into account the fact that the resources have to travel between the locations on a point-to-point basis. The limit problem as $N\to\infty$ is characterized and some applications to airfreight systems are shown.
\end{abstract}



\section{Introduction}\label{sint}

Locating a given number of points in a region, in order to fulfill a given optimization criterion, is a widely studied problem, and a large number of references on the field is available (see References), with many of them devoted to several applications to economy, urban planning, electronics, communication systems.

In the most common framework, a given bounded and closed region $\Omega\subset\R^d$ is considered, together with a given nonnegative function $\rho:\Omega\to\R^+$ which represents the distribution density of resources in $\Omega$. The goal is to concentrate the resources into a given number $N$ of points $x_1,\dots,x_N$ in an optimal way; assuming that the cost to move a unit mass from $x$ to $y$ is proportional to a suitable power $|x-y|^p$ of the distance, allows us to write the optimization problem as
\be\label{eloc}
\min\Big\{\int_\Omega\big(\dist(x,\Sigma)\big)^p\rho(x)\,dx\ :\ \Sigma\subset\Omega,\ \#\Sigma=N\Big\}.
\ee
Here $\Sigma$ is the unknown set of $N$ points to be determined, $\#\Sigma$ is the cardinality of $\Sigma$, and $\dist(x,\Sigma)$ is the distance function
$$\dist(x,\Sigma)=\min\big\{|x-y|\ :\ y\in\Sigma\big\}.$$

Problems of the form \eqref{eloc} are known as {\it location problems} and the existence of an optimal configuration is straightforward. On the contrary, in spite of its simplicity, the numerical computation of an optimal set $\Sigma$, when the number $N$ is large, presents big difficulties, essentially due to the fact that the cost in \eqref{eloc} admits a huge number of local minima, which prevents the use of fast gradient methods and makes necessary the implementation of global optimization methods that are in general much slower.

The asymptotic analysis, as $N\to+\infty$, has been performed (see for instance \cite{bjr02,bss11} and references therein) for problem \eqref{eloc} and gives important information about the limit density of optimal points $x_i\in\Sigma$. In Section \ref{sloc} we recall the main results about this issue.

The problem we deal with in the present paper is concerned with the optimal location of a given number $N$ of airports in a region $\Omega$. The airports collect the resources that are distributed in $\Omega$ with a known density $\rho(x)$; moreover, the goods travel between airports on a point-to-point basis, which provides an additional cost, called {\it routing cost}. The complete problem that comes out by adding location and routing costs will be discussed in Section \ref{spbm}. When the number $N$ of airports is large, we replace the location cost by its asymptotic counterpart and we discuss the corresponding first order necessary conditions of optimality. Finally, in Section \ref{snum} some numerical simulations are shown.

\section{The optimal location problem}\label{sloc}

The location problem consists in determining, given a region $\Omega\subset\R^d$ and a nonnegative function $\rho:\Omega\to\R$ describing the density of resources in $\Omega$, the position of a given number $N$ of points $x_1,\dots,x_N$ in $\Omega$ in order to minimize the work necessary to concentrate the resources in the points $x_i$. Assuming that the work to move a unit mass from $x$ to $y$ is proportional to $|x-y|^p$, the optimal location problem can be written as the minimization problem \eqref{eloc}.

When the number $N$ of points is large, as explained in the Introduction, the numerical computation of the optimal points $x_i$ is heavy. Therefore an asymptotic analysis as $N\to+\infty$, which provides, instead of the precise location of the points $x_i$, only their asymptotic density in $\Omega$, can be very helpful and with much lighter computational costs.

In the case under consideration $\Omega$ is a geographic region in $\R^2$ on which a density of resources $\rho$ is distributed; we assume that the density $\rho$ is known. We want to locate in $\Omega$ a given number $N$ of airports in the most efficient way according to a global cost that we are going to define. In this section we take into account only the location cost, while the routing cost will be considered in Section \ref{spbm}. The location cost consists in evaluating the work necessary to concentrate the resources distributed on $\Omega$ into the airports; we denote by $x_i$ the positions of the airports and we assume that the work to move a unit mass from a point $x$ to a point $y$ is proportional to $|x-y|^p$. 

Let us denote by $m_i$ the quantity of resources that will be concentrated at the point $x_i$ and by $\Omega_i$ the so-called {\it Voronoi cell} corresponding to $x_i$, that is the subregion of $\Omega$ that sends its resources to the point $x_i$. In other words, we have
$$m_i=\int_{\Omega_i}\rho(x)\,dx.$$
Then we have that the total cost to concentrate the resources spread on $\Omega_i$ into the airport $x_i$ is given by
$$A\int_{\Omega_i}|x-x_i|^p\rho(x)\,dx$$
where $A$ is a proportionality constant. Summing up over all the $N$ airports we have that the total cost is given by
$$A\sum_{i=1}^N\int_{\Omega_i}|x-x_i|^p\rho(x)\,dx$$
that can be also written in the form
\be\label{elocost}
A\int_\Omega\big(\dist(x,\Sigma)\big)^p\rho(x)\,dx
\ee
where $\Sigma$ is the unknown set of $N$ points to be determined. The most efficient choice of the positions of the airports, considering only the location cost, is then obtained by solving the minimization problem \eqref{eloc}.

When the number $N$ tends to $+\infty$, instead of looking at the precise positions $x_i$ in $\Omega$ of the airports, one will simply target at determining the limit density $\mu$ of the points $x_i$. In order to do it, we identify each set $\Sigma\subset\Omega$ of $N$ points with the measure
$$\mu_N=\frac{1}{N}\sum_{i=1}^N{\delta_{x_i}}.$$

If we assume that, up to a normalization, the density $\rho$ has a unitary total mass, the location cost \eqref{elocost} is proportional to the $p$-th power of the Wasserstein distance between the probabilities $\rho\,dx$ and $\mu_N$. The asymptotic analysis of the cost above has been performed (see for instance \cite{bjr02,bss11} and references therein) and we summarize here below the available results that, to be correctly stated, require the use of the $\Gamma$-convergence, a variational theory developed by De Giorgi and his school starting from the seventies.

\begin{figure}[ht]
\centerline{\includegraphics[height=5.0cm,width=10.0cm]{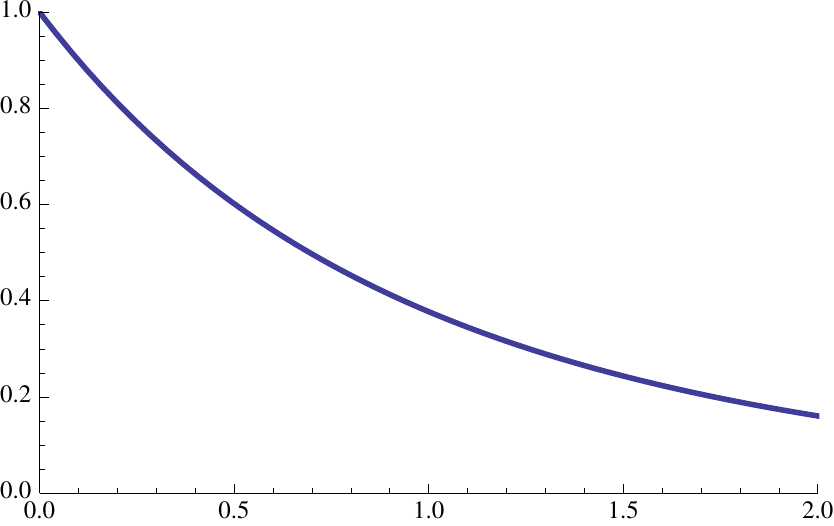}}
\caption{Plot of the value of $C_{p,2}$ for $p\in[0,2]$.}\label{fig1}
\end{figure}

When $N\to+\infty$ the cost \eqref{elocost} is asymptotically equivalent to the limit cost
\be\label{elimloc}
A\,C_{p,d}\,N^{-p/d}\,\int_\Omega\frac{\rho(x)}{\big(\mu(x)\big)^{p/d}}\,dx
\ee
expressed in terms of the limit density $\mu$ of points, where $C_{p,d}$ is a constant depending on the exponent $p$ and on the dimension $d$. It has to be noticed that in the integral above only the absolutely continuous part of $\mu$ has to be taken into account, neglecting the singular part. In the case $d=2$ the constant $C_{p,2}$ can be explicitly computed and we have
$$C_{p,2}=\int_E|x|^p\,dx$$
where $E$ is the regular hexagon of unitary area centered at the origin. For instance one has $C_{1,2}\sim0.377$ and $C_{2,2}\sim0.16$. A plot of the value of $C_{p,2}$ for $p\in[0,2]$ is given in Figure \ref{fig1}.

On the other hand, if we are interested not only in the location $x_i$ of the $i$-th airport but also in the mass $m_i$ that is there concentrated, instead of the measures $\mu_N$ above we have to consider the measures
$$\nu_N=\sum_{i=1}^Nm_i\delta_{x_i}$$
and the optimization problem is written in terms of the $p$-Wasserstein distance as
\be\label{ewassN}
\min\big\{W_p^p(\rho,\nu)\ :\ \#(\spt\nu)=N\big\}.
\ee
We notice that, without the normalization $\int\rho\,dx=1$, passing to the probability $\rho(x)\big/\int\rho\,dx$, the optimization problems \eqref{elimloc} and \eqref{ewassN} remain of the same form.

\section{Routing costs}\label{spbm}

In this section we consider some extra cost that occurs in the model of airports location. Once the resources, spread on $\Omega$ with density $\rho$, have been concentrated at the points $x_i$, with mass $m_i$ each, the goods will travel among airports that are connected on a point-to-point basis. We consider two situations: in the first one the routing cost depends linearly on the transported mass, while in the second the connection cost between $x_i$ and $x_j$ only depends on the distance $|x_i-x_j|$. In the following we denote by $V(x)$ the function $|x|^q$.

\subsection{Mass dependent routing costs}\label{smassc}

In this subsection we assume that the mass $m_i$ concentrated at the point $x_i$ is dispatched to the remaining points $x_j$ proportionally to the masses $m_j$; moreover, we assume that the cost to move a unit mass from a point $x$ to a point $y$ is proportional to $|x-y|^q$ for a suitable power $q$. Therefore, the cost to move the entire mass $m_i$ is
$$B\sum_jm_i\frac{m_j}{m}|x_i-x_j|^q$$
where $B$ is a proportionality constant and $m=\sum_jm_j=\int\rho\,dx$. Finally, the total routing cost is
$$\frac{B}{m}\sum_{i,j}m_im_j|x_i-x_j|^q.$$
If we write the routing cost in terms of the measure $\nu_N$ we obtain
$$\frac{B}{m}\int_\Omega\int_\Omega|x-y|^q\,d\nu_N(x)\,d\nu_N(y)=\frac{B}{m}\int_{\Omega\times\Omega}V(x-y)\,d(\nu_N\otimes\nu_N)$$
and the total cost taking into account location and routing terms gives the optimization problem
\be\label{etotm}
\min\Big\{AW_p^p(\rho,\nu)+\frac{B}{m}\int_{\Omega\times\Omega}V(x-y)\,d(\nu\otimes\nu)\ :\ \#(\spt\nu)=N\Big\}
\ee
The characterization of the limit problem as $N\to\infty$ in this case is easy and we can write it as
\be\label{etotmlim}
\min\Big\{AW_p^p(\rho,\nu)+\frac{B}{m}\int_{\Omega\times\Omega}V(x-y)\,d(\nu\otimes\nu)\Big\}
\ee
where the minimization above is intended in the class of all measures $\nu$ having the same total mass as $\rho$.

The necessary conditions of optimality for the optimization problem \eqref{etotmlim} can be obtained by differentiating the Wasserstein distance term (see \cite{busa09}) and the routing cost; we obtain
\be\label{ott}
A\phi+\frac{2B}{m}V*\nu=c\qquad\nu\hbox{-a.e.}
\ee
where $\phi$ is the Kantorovich potential for the transport from $\rho$ to $\nu$ and $c$ is a constant playing the role of the Lagrange multiplier of the mass constraint on $\nu$. In \eqref{ott} the measure $\nu$ appears in a very implicit way and can be determined only numerically. One connection between the Kantorovich pontential $\phi$ and the transport map $T$ from $\rho$ to $\nu$ is given by the Monge-Amp\`ere equation
$$\rho=\nu(T)\det(\nabla T).$$
Differentiating in $\eqref{ott}$ we obtain
$$A\nabla\phi+\frac{2B}{m}\nabla V*\nu=0$$
and $T(x)=x-\nabla\phi(x)$. Therefore we have the system
\be\label{sistema}
\begin{cases}
\ds A(x-T(x))+\frac{2B}{m}\nabla V*\nu=0\\
\rho=\nu(T)\det(\nabla T).
\end{cases}
\ee
In dimension 1 we can proceed by an iterative scheme, fixing an initial $\nu_0$ and obtaining $T_0$ from the first equation in \eqref{sistema}. Then we can recover $\nu_1$ by the second equation
$$\nu_1(T_0(x))=\frac{\rho(x)}{T_0'(x)}$$
and, assuming $T_0$ invertible,
$$\nu_1(y)=\frac{\rho(T_0^{-1}(y))}{T_0'(T_0^{-1}(y))}.$$
We can now proceed by iterating the scheme above.

\begin{exam}
In this particular example we can find an explicit solution taking $n=1$, $p=2$, and $V(s)=|s|^2$. If we suppose that the barycenter of $\nu$ is in the origin, we obtain:
$$V*\nu=mx^2+\int y^2\,d\nu(y)$$
so that
$$A\phi'(x)+4Bx=0$$
which gives
$$\phi'(x)=-\frac{4B}{A}x\qquad\hbox{and}\qquad T(x)=\Big(1+\frac{4B}{A}\Big)x.$$
Putting this in the $1$-dimensional Monge-Amp\`ere equation and indicating by $v$ the density of $\nu$, we obtain
$$\rho(x)=v((1+4B/A)x)\Big(1+\frac{4B}{A}\Big),$$
and changing variables,
$$v(y)=\frac{1}{1+4B/A}\rho\Big(\frac{y}{1+4B/A}\Big).$$
\end{exam}

\subsection{Mass independent routing costs}\label{snomassc}

In this subsection we assume that the cost to connect the airports located at the points $x_i$ and $x_j$ does not depend on the transported mass and amounts simply to $K|x_i-x_j|^q$ where now the constant $K$ is the cost of flying along a unit distance. In this case it is more convenient to use the probability measures $\mu_N$ introduced in Section \ref{sloc} which provide the routing cost in the form
$$K\sum_{i,j}|x_i-x_j|^q=KN^2\int_{\Omega\times\Omega}V(x-y)\,d(\mu_N\otimes\mu_N).$$
Taking into account the asymptotic expression of the location cost given in \eqref{elimloc}, we obtain the optimization problem
$$\min\Big\{A\,C_{p,d}\,N^{-p/d}\int_\Omega\frac{\rho(x)}{\big(\mu(x)\big)^{p/d}}\,dx+KN^2\int_{\Omega\times\Omega}V(x-y)\,d(\mu\otimes\mu)\Big\}$$
where now $\mu$ runs in the class of all probabilities on $\Omega$. Setting $\eps=A\,C_{p,d}\,N^{-2-p/d}/K$ we are now faced with the problem
\be\label{etotnom}
\min\Big\{F_\eps(\mu):=\eps\int_\Omega\frac{\rho(x)}{\big(\mu(x)\big)^{p/d}}\,dx+\int_{\Omega\times\Omega}V(x-y)\,d(\mu\otimes\mu)\Big\}.
\ee

The necessary conditions of optimality for problem \eqref{etotnom} simply follow by differentiation of the cost functional and give:
\be\label{econveq}
\eps\rho\,\frac{p}{d}\,\mu^{-1-p/d}+2V*\mu=c
\ee
where $*$ denotes the convolution operator and $c$ is a constant coming from the mass constraint on $\mu$.

When $\eps\to0$ the optimal densities $\mu_\eps$ of problem \eqref{etotnom} tend to a Dirac mass $\delta_{x_0}$ for a suitable point $x_0$. In order to identify the limit problem as $\eps\to0$, and so to identify the point $x_0$ around which the optimal densities $\mu_\eps$ concentrate (it can be seen as the {\it main hub} of the airports system), it is convenient to rescale the cost above dividing it by its minimum value. Considering the measures
$$\mu=\delta\frac{1}{|\Omega|}+(1-\delta)\frac{1_{B_r(x_0)}}{|B_r(x_0)|}\qquad\hbox{with }r^q\ll\delta$$
a simple calculation provides for the minimal cost of problem \eqref{etotnom}
$$\min F_\eps\sim C\eps\delta^{-p/d}+\delta$$
for a suitable constant $C$. Optimizing with respect to $\delta$ the quantity above we obtain $\delta\sim\eps^{1/(1+p/d)}$ so that
$$\min F_\eps\sim\eps^{1/(1+p/d)}.$$
Therefore the rescaled functionals become
$$G_\eps(\mu)=\eps^{(p/d)/(1+p/d)}\int_\Omega\frac{\rho(x)}{\big(\mu(x)\big)^{p/d}}\,dx+\eps^{-1/(1+p/d)}\int_{\Omega\times\Omega}V(x-y)\,d(\mu\otimes\mu).$$
Note that the optimal measures for $F_\eps$ and for $G_\eps$ are the same. 

In order to characterize the asymptotic behavior of the minimizing sequences $(\mu_\eps)$ we will compute in the next section the $\Gamma$-limit of the sequence of functionals $(G_\eps)$. The general theory of $\Gamma$-convergence (see for instance \cite{dm93}) then provides the identification of the main hub $x_0$ around which the measures $\mu_\eps$ tend to concentrate.

\section{The $\Gamma$-convergence result}\label{sgam}

First of all we notice that, due to the presence of the coefficient $\eps^{-1/(1+p/d)}$ in front of the routing term, the $\Gamma$-limit on a measure $\mu$ will be $+\infty$ whenever $\int_{\Omega\times\Omega}V(x-y)\,d(\mu\otimes\mu)\ne0$. Therefore, we may limit ourselves to analyze only the measures for which the routing term vanishes, i.e. the Dirac masses $\mu=\delta_{x_0}$.

It is convenient to set
$$\alpha=\frac{p/d}{1+p/d}\;,\qquad\beta=\frac{1}{1+p/d}\;;$$
notice that $\alpha+\beta=1$ and that $\alpha=\beta p/d$. We will show that the $\Gamma$-limit of the sequence of functionals $G_\eps$, computed on the dirac mass $\delta_{x_0}$ and with respect to the weak* convergence of measures, coincides with the functional
$$H(\delta_{x_0})=A\int_\Omega\big(\rho(x)\big)^\beta|x-x_0|^{\alpha q}\,dx\qquad\hbox{where }A=\Big(1+\frac{p}{d}\Big)\Big(\frac{2d}{p}\Big)^\alpha.$$

\subsection{The $\Gamma$-limsup inequality}\label{glimsup}

In order to obtain a $\Gamma$-limsup inequality we have to choose a suitable sequence $\mu_\eps\rightharpoonup\delta_{x_0}$ and compute the limit of $G_\eps(\mu_\eps)$. We take
$$\mu_\eps=\eps^\beta\phi+\Big(1-\eps^\beta\int_\Omega\phi\,dx\Big)\delta_{x_0}$$
where the function $\phi$ will be chosen later. Then $\mu_\eps$ is a probability measure and we have
$$\begin{array}{lll}
G_\eps(\mu_\eps)&\ds=\eps^\alpha\int_\Omega\frac{\rho(x)}{\eps^{\beta p/d}\phi^{p/d}}\,dx
+\eps^{-\beta}\int_{\Omega\times\Omega}\eps^{2\beta}V(x-y)\phi(x)\phi(y)\,dxdy\\
&\ds\hskip3truecm+\eps^{-\beta}\int_\Omega2\eps^\beta\Big(1-\eps^\beta\int_\Omega\phi\,dx\Big)V(x-x_0)\phi(x)\,dx\\
&\ds=\int_\Omega\Big[\frac{\rho(x)}{\phi^{p/d}}+2V(x-x_0)\phi\Big]\,dx
+\eps^\beta\int_{\Omega\times\Omega}V(x-y)\phi(x)\phi(y)\,dxdy\\
&\ds\hskip3truecm-2\eps^\beta\int_\Omega\phi\,dx\int_\Omega V(x-x_0)\phi(x)\,dx
\end{array}$$
which gives
$$\lim_{\eps\to0}G_\eps(\mu_\eps)=\int_\Omega\Big[\frac{\rho(x)}{\phi^{p/d}}+2V(x-x_0)\phi\Big]\,dx.$$
We choose now $\phi$ in order to minimize the quantity at the right-hand side. An easy computation gives
$$\phi(x)=\Big(\frac{p}{2d}\frac{\rho(x)}{V(x-x_0)}\Big)^\beta$$
which implies
$$\lim_{\eps\to0}G_\eps(\mu_\eps)=H(\delta_{x_0}).$$

\subsection{The $\Gamma$-liminf inequality}\label{gliminf}

In order to conclude that the $\Gamma$-limit of the functionals $G_\eps$ is the functional $H$ it remains to show the $\Gamma$-liminf inequality, which amounts to prove that for every sequence $\mu_\eps\rightharpoonup\delta_{x_0}$ we have
\be\label{infineq}
\liminf_{\eps\to0}G_\eps(\mu_\eps)\ge H(\delta_{x_0}).
\ee
The following lemma will be useful.

\begin{lemm}\label{singm}
Let $\mu$ be a measure on $\Omega$ that is singular with respect to the Lebesgue measure and let $\mu_n\rightharpoonup\mu$. Then there exists a sequence of open sets $(A_n)$ such that:
\begin{itemize}
\item[i)]$|A_n|\to0$;
\item[ii)]$\mu_n(\Omega\setminus A_n)\to0$ (hence $\mu_n\lfloor A_n\rightharpoonup\mu$).
\end{itemize}
\end{lemm}

\begin{proof}
Since $\mu$ is singular, it is concentrated on a measurable set $S$ with $|S|=0$ and there is a sequence of open sets $(A_k)$ containing $S$ and such that $|A_k|\to0$. Since $\Omega\setminus A_k$ are closed sets, we have for every $k$
$$\limsup_{n\to\infty}\mu_n(\Omega\setminus A_k)\le\mu(\Omega\setminus A_k)=0.$$
A diagonal argument achieves the proof.
\end{proof}

Take now a generic sequence $\mu_\eps\rightharpoonup\delta_{x_0}$, denote by $u_\eps(x)$ the density of the absolutely continuous part of $\mu_\eps$ with respect to the Lebesgue measure, and let $A_\eps$ be the open sets provided by Lemma \ref{singm}. Define
$$\mu^1_\eps=\mu_\eps\lfloor A_\eps,\qquad\mu^2_\eps=\mu_\eps\lfloor A^c_\eps.$$
We have
$$\begin{array}{lll}
G_\eps(\mu_\eps)&\ds=\eps^\alpha\int_\Omega\frac{\rho}{u_\eps^{p/d}}\,dx
+\eps^{-\beta}\Big[\int_{\Omega\times\Omega}V(x-y)\,d(\mu^1_\eps\otimes\mu^1_\eps)+\int_{\Omega\times\Omega}V(x-y)\,d(\mu^2_\eps\otimes\mu^2_\eps)\\
&\hskip3truecm\ds+\int_{\Omega\times\Omega}2V(x-y)\,d(\mu^1_\eps\otimes\mu^2_\eps)\Big]\\
&\ds\ge\int_\Omega\Big[\eps^\alpha\frac{\rho}{u_\eps^{p/d}}+\eps^{-\beta}2(V*\mu^1_\eps)u_\eps1_{A^c_\eps}\Big]\,dx
\end{array}$$
where we used the fact that $\int_{\Omega\times\Omega}V(x-y)\,d(\nu\otimes\nu)\ge0$ for every measure $\nu$ and that $\mu^2_\eps\ge u_\eps1_{A^c_\eps}\,dx$. Using the Young inequality
$$X\eps^\alpha+Y\eps^{-\beta}\ge\frac{X^\beta Y^\alpha}{\alpha^\alpha\beta^\beta}\qquad\hbox{for }\alpha+\beta=1,$$
we obtain
$$\begin{array}{ll}
G_\eps(\mu)&\ds\ge\int_\Omega\frac{1}{\alpha^\alpha\beta^\beta}\Big(\frac{\rho}{u_\eps^{p/d}}\Big)^\beta\big(2(V*\mu^1_\eps)u_\eps1_{A^c_\eps}\big)^\alpha\,dx\\
&\ds=A\int_{A^c_\eps}\rho^\beta(V*\mu^1_\eps)^\alpha\,dx.\\
\end{array}$$
Since $|A_\eps|\to0$ and $(V*\mu^1_\eps)(x)\to V(x-x_0)$ we finally obtain \eqref{infineq}.

As a conclusion, the $\Gamma$-limit computation is achieved and the optimal {\it main hub} for the limit location-routing problem of \eqref{etotnom} is located at the point $x_0$ which minimizes the quantity

\begin{equation}\label{ex0} 
\int_\Omega\big(\rho(x)\big)^{1/(1+p/d)}|x-x_0|^{q(p/d)/(1+p/d)}\,dx.
\end{equation}
Note that this minimization problem for $x_0$ is of the form of a Torricelli optimal location problem with suitable exponents.

\section{Some numerical simulations}\label{snum}

In this section we perform some numerical simulations, based on the results of the previous section, that can be applied to real cases. In the first subsection, some 1-dimension and 2-dimension examples will be presented for different routing costs and density functions and varying $\eps$.
The last subsection introduces an application of the model to the USA airfreight system in order to compare the results with the current location of US airfreight hubs.

We are interested to find minimizer of the functional representing the sum of location and routing costs and the $x_0$ minimizing the $\Gamma$-limit functional $H(\delta_{x_0})$. In fact, one of the properties of $\Gamma$-convergence is the convergence of minima, so if $H$ is the $\Gamma$-limit of $G_\eps$, the limit of minimizers of $G_\eps$ is a minimizer of $H$. We will find numerically the optimal $\mu_\eps$ and we observe that for small $\eps$ they are close to a Dirac mass at a suitable point $x_0$, according to the previous section's result.

In Section \ref{spbm} we found an optimality condition for the minimizers of $F_\eps$. Unfortunately, condition \eqref{econveq} does not admit an explicit solution, so we approximate it numerically. More specifically, we approximate the minimizer via an iterative scheme: we start from the uniform distribution with total mass $1$ and then we define the iteration term according to the necessary condition:
\be\label{itterm}
\begin{cases}\mu_0=\mathcal{U}(\Omega)\\
\mu_{n+1}=\Big(\frac{\eps\rho}{c+V*\mu_n}\Big)^{p/d +1}.\end{cases}
\ee
Here the Lagrange multiplier $c$, according to condition \eqref{econveq}, is proportional to $\eps^{1/(1+p/d)}$.

\subsection{1-D and 2-D examples}\label{examples}

\begin{figure} [ht]
\centering
\includegraphics[height=9cm]{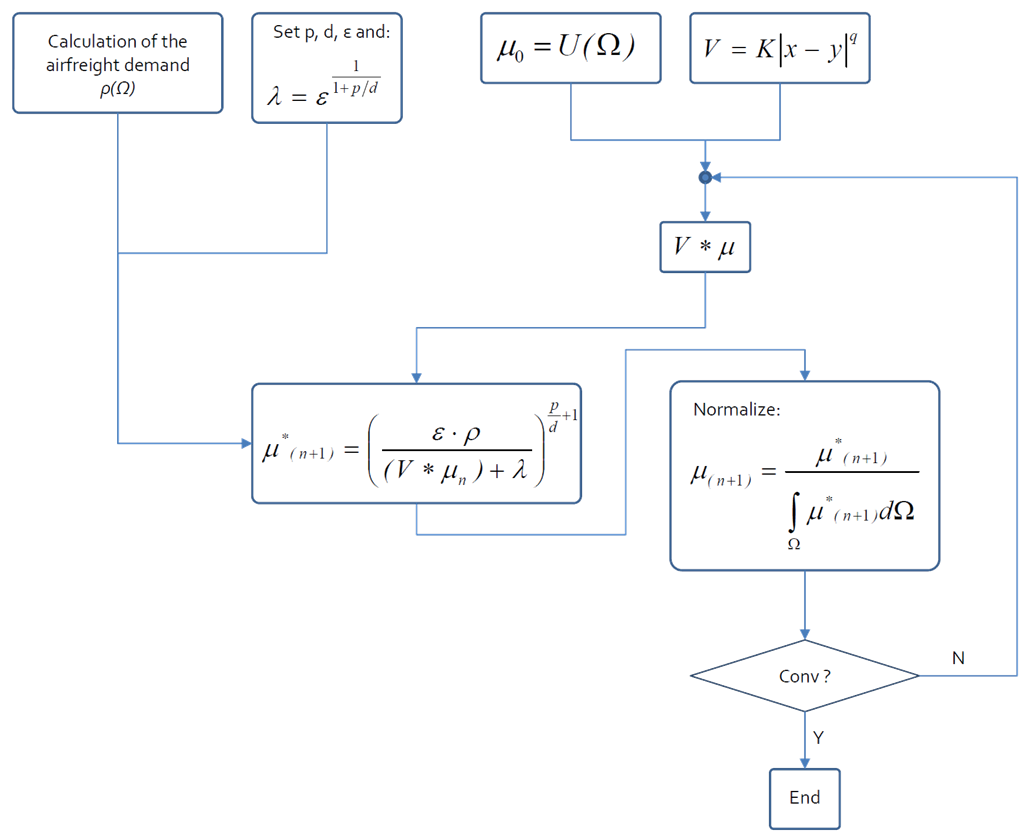}
\caption{numeric procedure for the determination of probability distribution $\mu$}
\label{proc}
\end{figure}

In the one-dimensional case, the domain $\Omega$ is the interval $[-1,+1]$ that is discretized in order to solve numerically the problem. Consequently, both the functions $\rho$ and $\mu$ are expressed through an array of values in correspondence of the discretization points. At each step of the convergence procedure shown in Figure \ref{proc}, $\mu_{n+1}$ is obtained from the relationship \eqref{itterm} and then normalized to a probability measure.

\begin{figure}
\centering
\includegraphics[height=8cm, width=10.5cm]{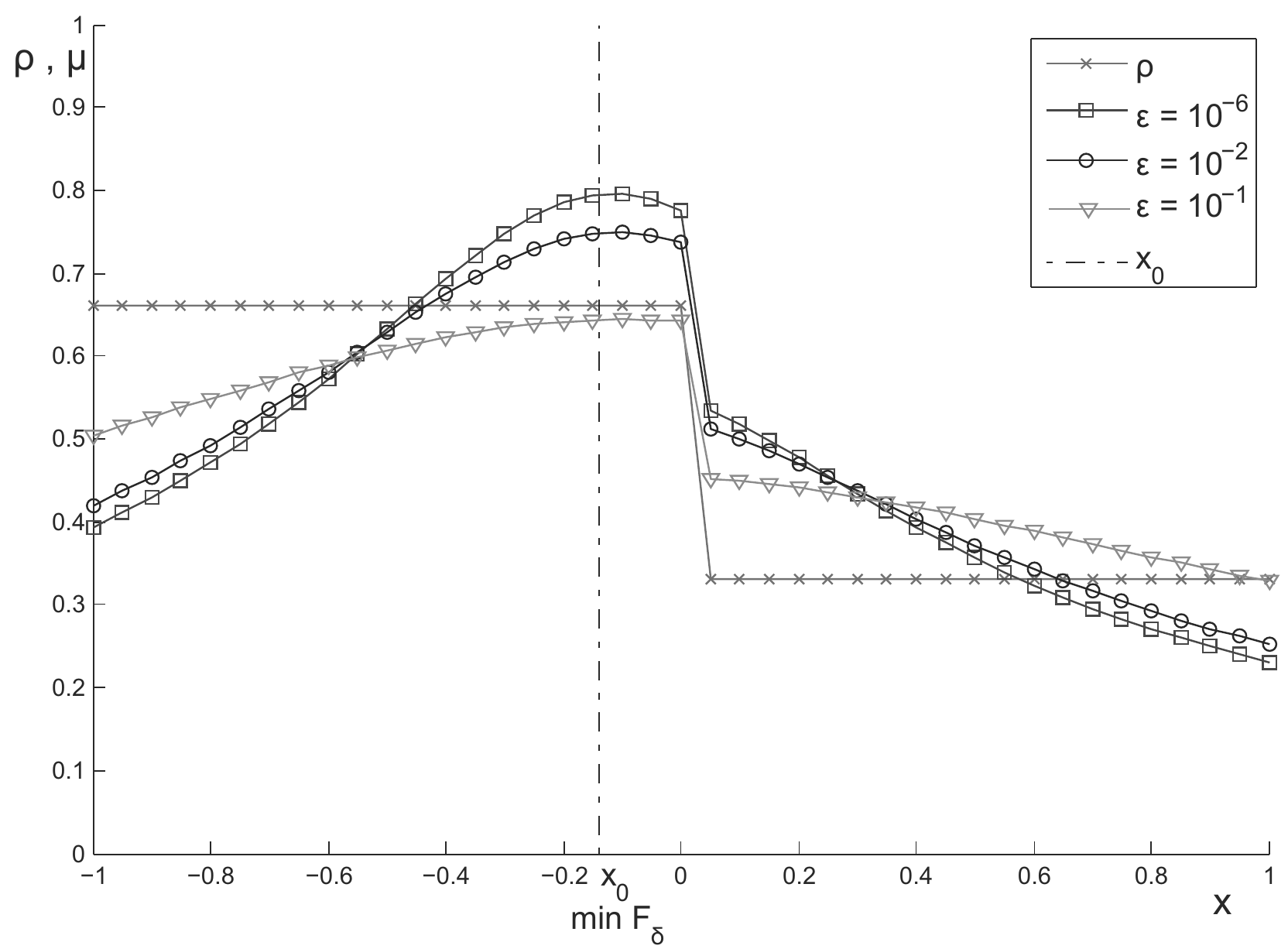}
\caption{results of the first simulation (asymmetric population)}
\label{resultasym}
\end{figure}

The first simulation, reported in Figure \ref{resultasym}, is related to a non-symmetric distribution of population density $\rho$ and a quadratic routing cost function:
$$\rho(x)=\begin{cases}
1&\hbox{if }x\in[-1,0]\\
2&\hbox{if }x\in[0,+1],
\end{cases}
\qquad V=|x-y|^2,\qquad p=1.$$

We assume that the convergence is reached when the maximum error between the values of $\mu_n$ and $\mu_{n+1}$ is less than 2\%. In this conditions, about 10 iterations are requested to solve the problem and the computational time results to be proportional to the number of point used to discretize the domain; when it is divided into 200 steps, the calculation time is about 100 sec.

The results for different values of $\eps$ coefficient show that, as the $\eps$ decreases, the limit density $\mu$ tends to have a concentration centred on a single point. At the limit as $\eps\to0$, the density $\mu$ becomes a Dirac mass located at the point $x_0$ that minimizes the functional \eqref{ex0}. The convergence towards the limit conditions of $\eps\to0$ is slow and it cannot be reached numerically because the onset of numerical errors below the value of $\eps\simeq10^{-4}$. Therefore, the routine can be completed by calculating also the value of the functional $H(\delta_{x_0})$ reported in (\ref{ex0}), and find the point $x_0$ of minimum. In this case, the minimum of the functional \eqref{ex0} can be found explicitly:
\be
\begin{split}
H(\delta_{x_0})&=\int_{-1}^1\sqrt{\rho(x)}|x-x_0|\,dx\\
&=\sqrt{2}\int_{-1}^{x_0}(x_0-x)\,dx+\sqrt{2}\int_{x_0}^0(x-x_0)\,dx+\int_0^1(x-x_0)\,dx\\
&=\sqrt2\,x_0^2+x_0(\sqrt2-1)+\frac{1}{2}(\sqrt2+1)
\end{split}
\ee
which gives
$$x_0=\frac{\sqrt2-2}{4}\simeq-0.146.$$

The analytical solution equals to the value determined by the numerical procedure that is also reported in Figure \ref{resultasym}.

The population $\rho$ often can have an uneven distribution among the domain, and therefore an adequate function is requested in order to model correctly this aspect. A first solution can be provided by treating this distribution as a sum of M Gaussian functions:
\be\label{peakfun}
\rho(x)=\sum_{j=1}^M A_j e^{-B_j|X_j-x|^2}
\ee 
where the coefficients $A_j$,$B_j$,$X_j$ are used to set respectively the height, the width and the position of the $j$-th peak.

\begin{figure}
\centering
\includegraphics[height=7.5cm, width=10.5cm]{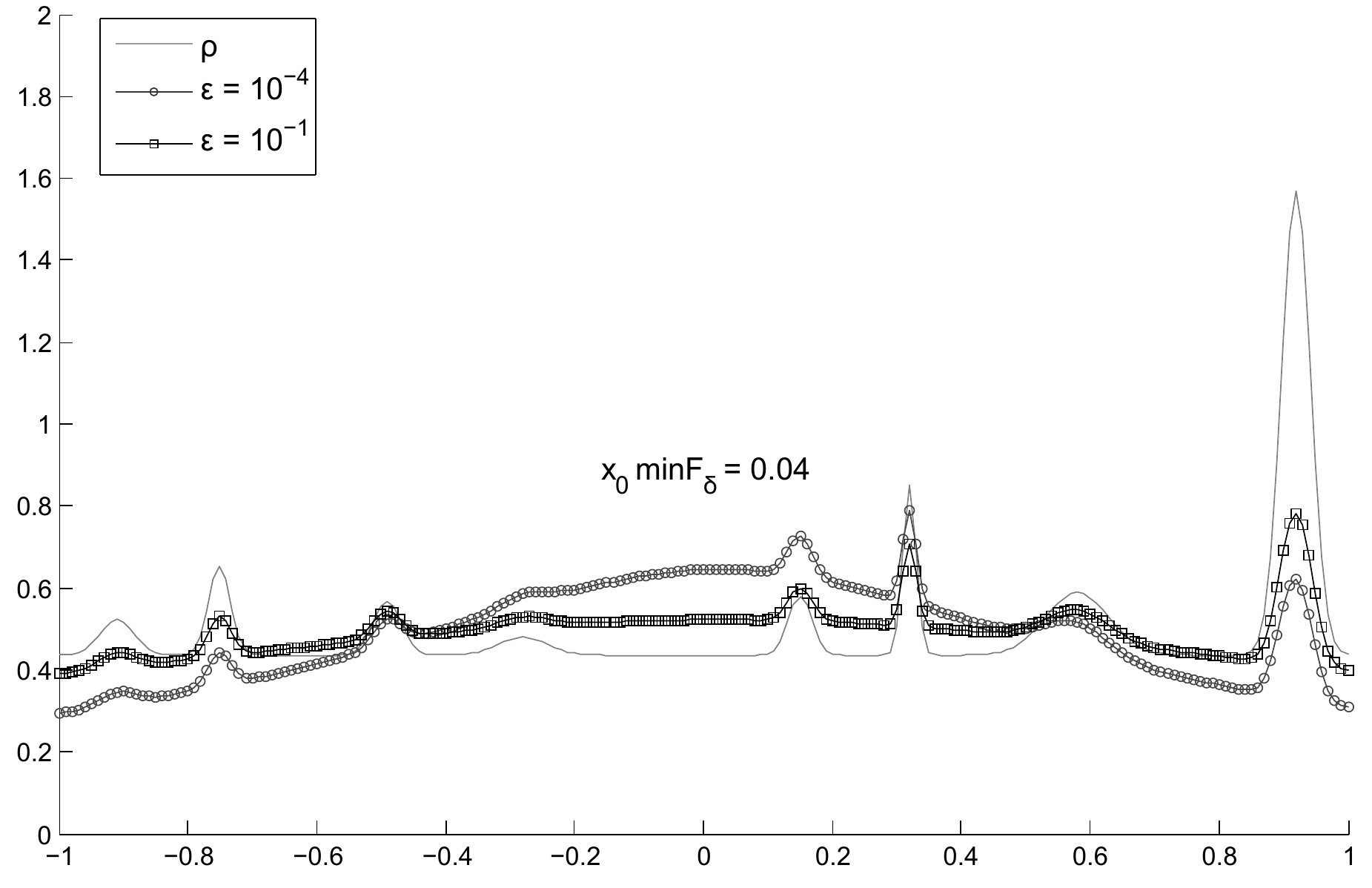}
\caption{second simulation: population modelled through a sum of Gaussian function}
\label{1Dpeak}
\end{figure}

The results reported in Figure \ref{1Dpeak} refer to the case of a population $\rho$ with 8 peaks of different position, height, and area of influence (width); also in this case the simulations have been conducted with two different values of the coefficient $\eps$.

The density of probability $\mu$ is highly dependent by the values of the coefficients; the point $x_0$ is in this case numerically determined: 
$$x_0\simeq+0.044.$$

As the $\eps$ decrease, the influence of the routing cost becomes larger on despite of the location ones so that the system tends to minimize the airport distance. The large differences between the solutions remarks the importance in choosing a value of the coefficient as realistic as possible.
Moreover, one can note that the computational time is not affected by the ``complication level'' of the population function but only by the used discretization step. The effective decisional process related to the facilities (the airports) location, can be done in a post-processing phase: in this way, we can decide how many airports can be located in a given region, proportionally to the area limited by the density distribution; for example, the numbers on the $X$-axis of Figure \ref{1Dpeak} equals to the airport on each step (each step length is 0.2).

\begin{figure} [h]
\centering
\includegraphics[height=7 cm]{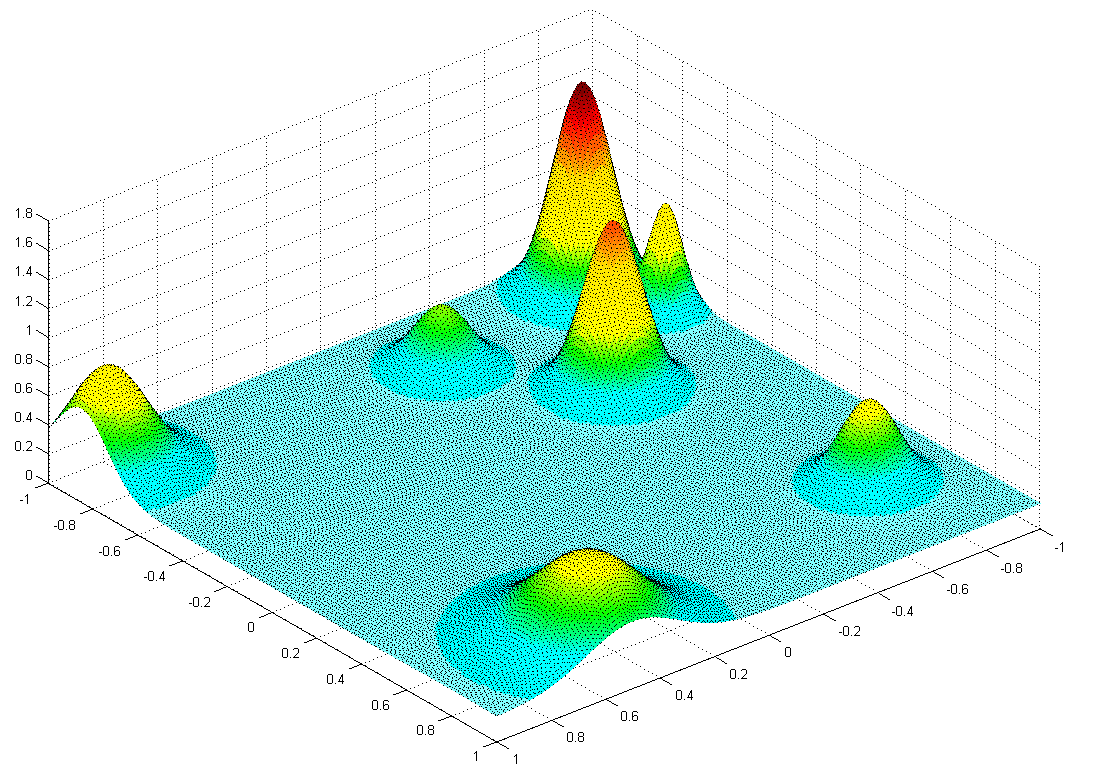}
\caption{2-D simulation: population $\rho$}
\label{2Drho}
\end{figure}

The routine has been applied also in the 2-D case, considering a correspondent peaks distribution shown in Figure \ref{2Drho}. We remember the Gaussian function in the case of two variables and also both the routing cost function and the exponent of location:
$$\rho(x,y)=\sum_{j=1}^M A_j e^{-B_j(|X_j-x|^2+|Y_j-y|^2)},\qquad V(x)=|x|^{0.5},\qquad p=1.$$

Although the calculation procedure does not change, the computational time results much higher than in the previous cases because of the great number of points requested to discretize properly the domain. Neverthless, it remains notably lower than the ones of the common Operating Research models (about 2200 sec. when the domain is divided into 1600 cells).

\begin{figure} 
\centering
\includegraphics[height=5.5 cm]{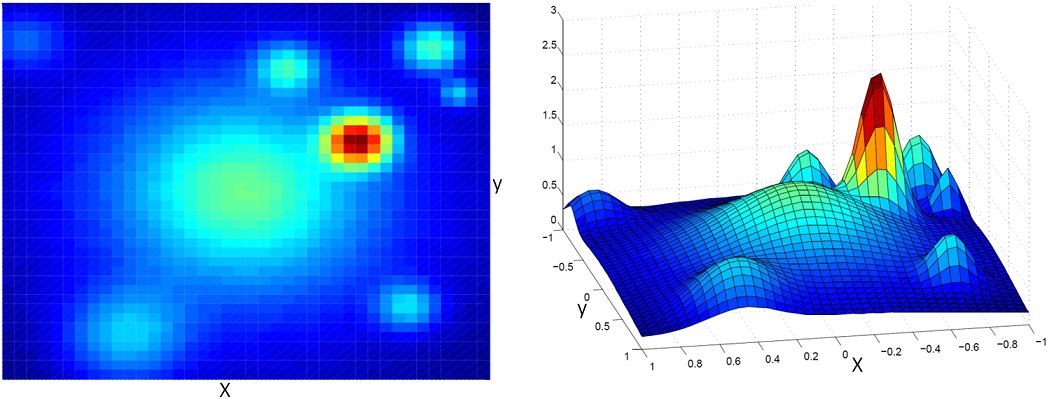}
\caption{2-D simulation: result}
\label{2Dresult}
\end{figure}

The result shows that the probability density follows the shape of the initial population $\rho(x,y)$ (we can note that the exponent q=0.5 determines the minor importance of the routing costs on despite of location ones) the point of maxima can be observed near the central peaks where the effect of both the location and routing terms are summed.

\begin{figure} [h]
\centering
\includegraphics[height=6 cm]{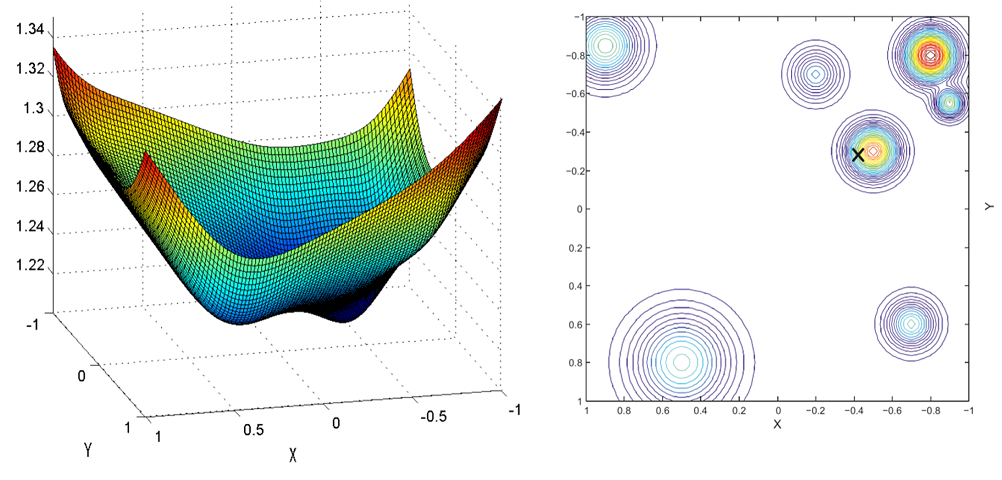}
\caption{2-D main hub result}
\label{2Dfunc}
\end{figure}

When the population distribution becomes irregular the position $x_0$ of the ``main hub'' cannot be estimated immediately but the functional \eqref{ex0} can be easily computed and its minimum can be founded. For the same 2-D case, the values of the functional are shown in Figure \ref{2Dfunc} while its minimum point is depicted in Figure \ref{2Dfunc} together with the level curves of the population $\rho$.

\subsection{Application to the US airfreight system}\label{application}

Two main problems will be faced in order to apply location-routing models to real cases:

\begin{itemize}
\item Location and Routing terms are related with ground and air transportation costs respectively. Preliminarily, we can suppose a linear dependence of ground cost with the transport distance but the same assumption becomes not valid in the case of air transportation.

\item The distribution of population $\rho$ identifies the airfreight demand among the domain; data are directly available only for some areas (occidental countries) while in most cases an extrapolation from some socio-economic data is needed.
\end{itemize}

Therefore, the aim of the present section is to set the coefficients and exponents appearing in \eqref{etotnom} in such the way the terms of cost functional reflect as realistic as possible the dynamics of the real world. The numerical routine will be finally applied on the US domain that will be represented as a polygon on a Cartesian system.

\subsubsection{Routing and airfreight cost}\label{routair}
\begin{figure} [h]
\centering
\includegraphics[height=5 cm]{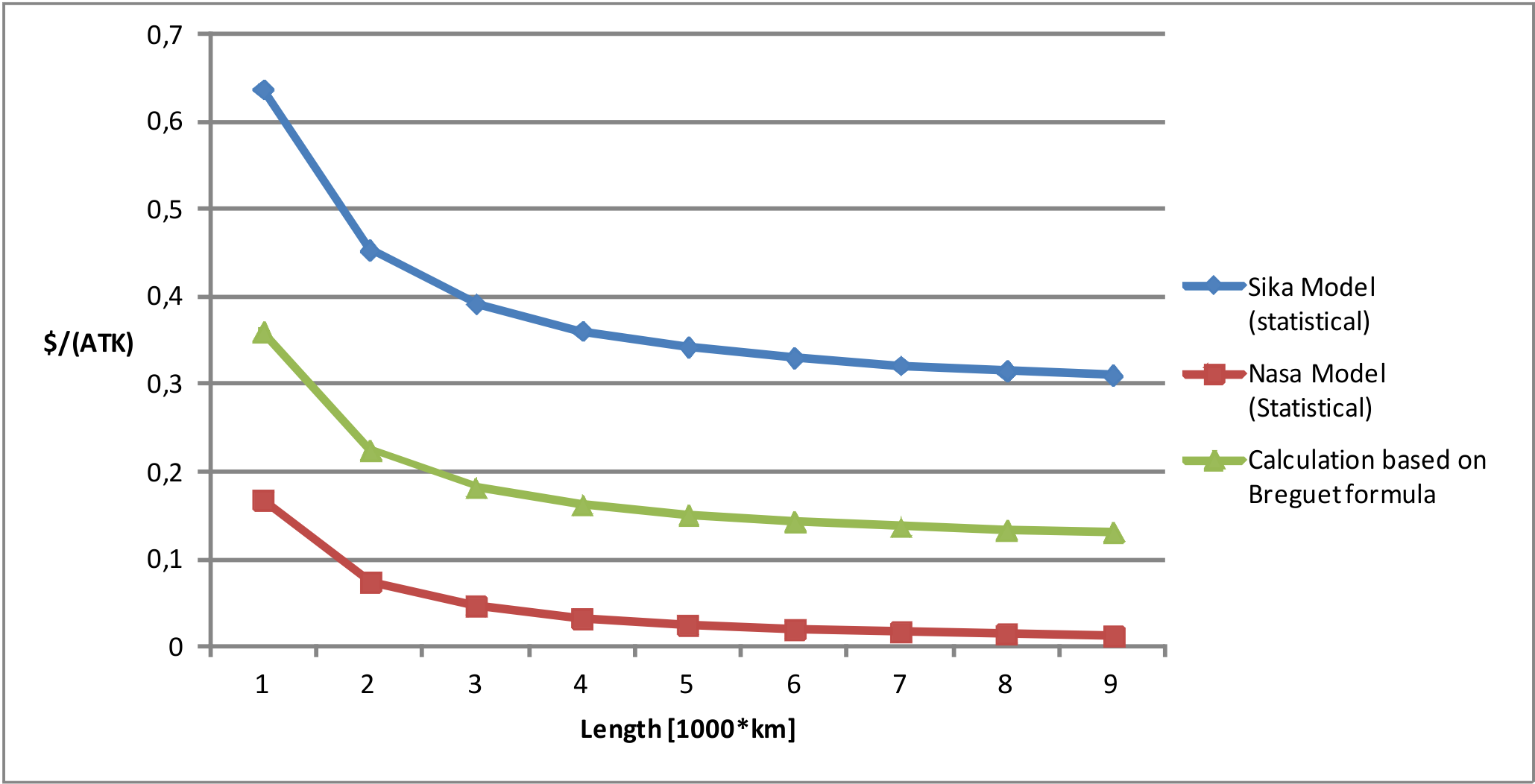}
\caption{Variation of air cost with transport distance}
\label{c_comp}
\end{figure}

In the functional \eqref{ex0}, the routing costs can be modelled as a function of the transport distance, through the general power relationship $V(x)=K|x|^q$ by simply setting the coefficient $K$ and the exponent $q$. In the case of airfreight, most of cost terms depend strongly by the \emph{economies of scale} in which the company operates (countries connected, commercial accordances, kind of service done, aircraft used) so that it often is not possible to determine a general function that could be valid in every case. Nevertheless, if we cannot determine an explicit function, we can determine its "shape" by supposing that the part of cost variable with the transport distance, is mainly related to the fuel consumption during flight. In a first approximation, the amount of fuel required for a given mission can be determined (for example in the case of constant power aircraft) by using the so-called Breguet relations: 

\be \label{eqbreguet} 
Cost_{fuel}\propto W_{fuel}= 1-e^{\frac{-Range *k_c}{\eta_p E}},
\ee
where the $Range$ equals to the transport distance, $E$ is the aerodynamic efficiency of the considered aircraft, and $k_c$ and $\eta_p$ are respectively the Specific Fuel Consumption and the propeller efficiency: since these parameters are all known for each aircraft and engine, the amount of fuel and its cost can be calculated in dependence of the Range flown.
The operating costs are usually reported in terms of \emph{Costs per Unit of carried mass and flown distance, $Cost/(Ton \cdot Km)$}(simply by dividing by the total payload and the transport distance) and the results of this procedure deriving from Breguet, has been compared in Figure \ref{c_comp} with some statistical models (\cite{nasacr} and \cite{sika}) that use a regression of both historical data about existing freighter and data collection of the financial report of transport companies.

The Curves in Figure \ref{c_comp} have a similar shape and they differ only by a translating coefficient that can be related to the different economies of scale which data are extrapolated from.
Moreover, the $Cost/Ton$ can be determined by integrating the curves in Figure \ref{c_comp} so that finally a suitable value of the exponent $q$ is determined:
$$q=0.7$$

The value results lower than $1$ because, conceptually, the air cannot be considered as a constant mass transport: in the case of existing aircraft in fact the Weight of embarked fuel represents until 30-40\% of the total mass so that the flight condition and, consequently, also the fuel vary notably during cruise (as the aircraft is lightening, the burn fuel decreases).

\subsubsection{Modeling the airfreight demand}\label{roudemand}

The airfreight demand is often not directly measurable so that also the initial population $\rho$ has to be properly modelled. The identification of the socio-economic parameters affecting the airfreight demand, results very difficult; also in this case the models are highly affected by the economies of scale and the geographical region on which the air transport is operated.
In the present study, we assume that the airfreight depends by some socio-economic parameters in a likely linear regression as proposed in \cite{alkaabi}. In the model proposed, the airfreigth demand in some point of the domain is obtained by the following expression:
\be\label{demandfun}
ln(AF)=C_0 + C_1 PC +C_2 TSE + C_3 TSL +C_4MD+C_5 HT
\ee
where 
\begin{itemize}
\item $C_0,..,C_5$: coefficient coming from linear regression of economic data
\item $AF$: volume of airfreight demand ($TON$)
\item $PC$: per capita personal income (\$1,000)
\item $TSE$: traffic shadow effect. In first approximation, this parameter will not considered in order to avoid any iterative process also for the input data.
\item $TSL$: transportation-shipping-logistics employment market share (\%)
\item $MD$:$\#$ of medical diagnostic establishments
\item $HT$: average high-tech employee wage (\$1,000)
\end{itemize}

\begin{figure} [h]
\centering
\includegraphics[height=5.7 cm]{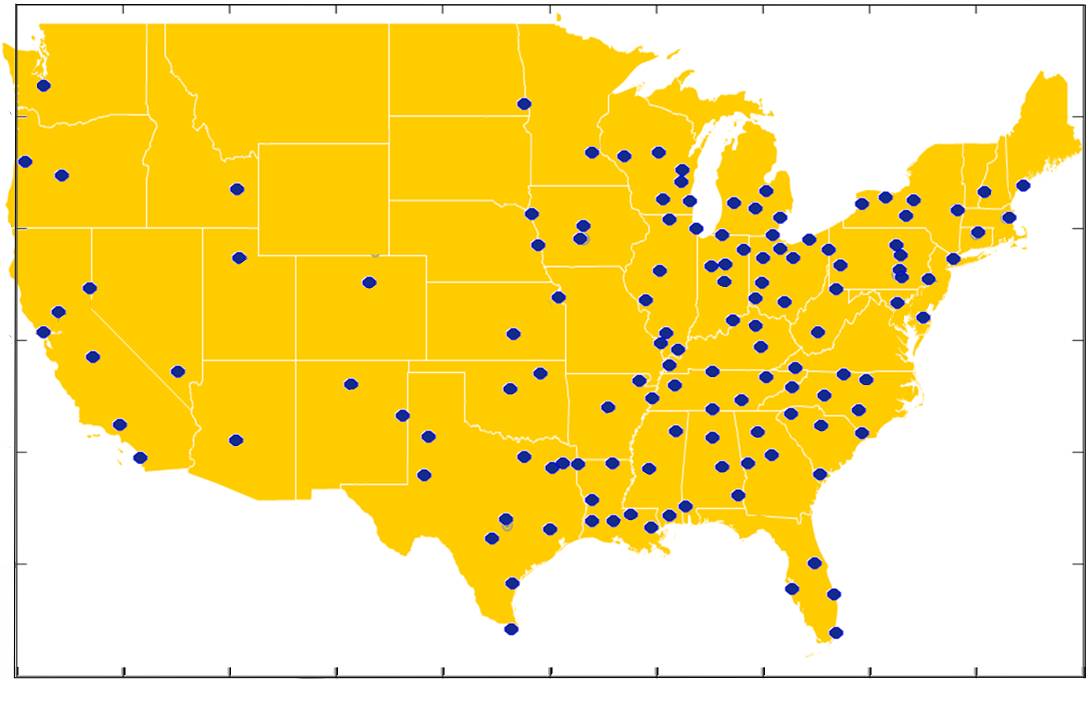}
\caption{centroids of MSAs among the US.}
\label{map}
\end{figure}

An airport has a relatively small catchment region (cities or districts), so that on the "ground side", the airfreight demand has influence on a very small area on despite of the worldwide dimension of the air transport; for this reason the airfreight function must to be refined also if the domain is very large and, consequently, the required discretization step is small.
In this contest,the socio-economic data in \eqref{demandfun} are extrapolated from common statistical reports (National Bureau for USA or Eurostat for E.U.) for each metropolitan districts, so that a detailed function can be easily determined.

The points used to define the spatial distribution of airfreight demand, are reported in Figure \ref{map} in blue dots, and their position refers to the centroids of the so called \emph{metropolitan statistical area} of the US territory.

\begin{figure} [h]
\centering
\includegraphics[height=5.5 cm]{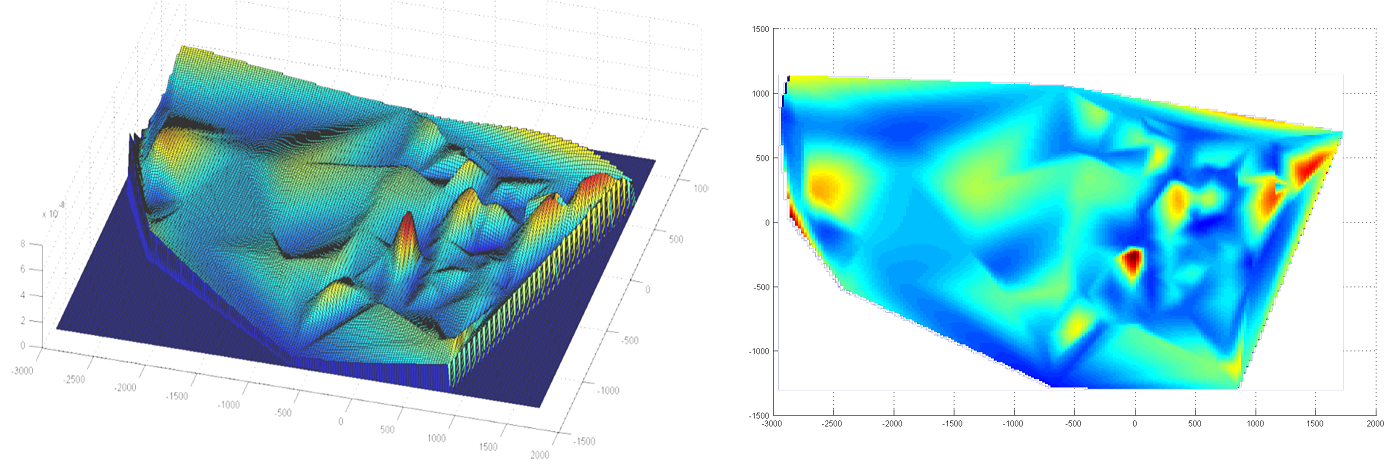}
\caption{model of the airfreight demand in Us}
\label{3ddemand}
\end{figure}

Since data are know in correspondence of these points, the $\rho$ function is then obtained through a cubic interpolation with a matlab routine. In Figure \ref{3ddemand} is reported the $\rho$ obtained by this procedure: it has an uneven distribution, peaks are concentrated in very rich or very populated regions and their area of influence is relatively small.

\subsubsection{Results}\label{ssresult}

\begin{figure} [!]
\centering
\includegraphics[height=8cm]{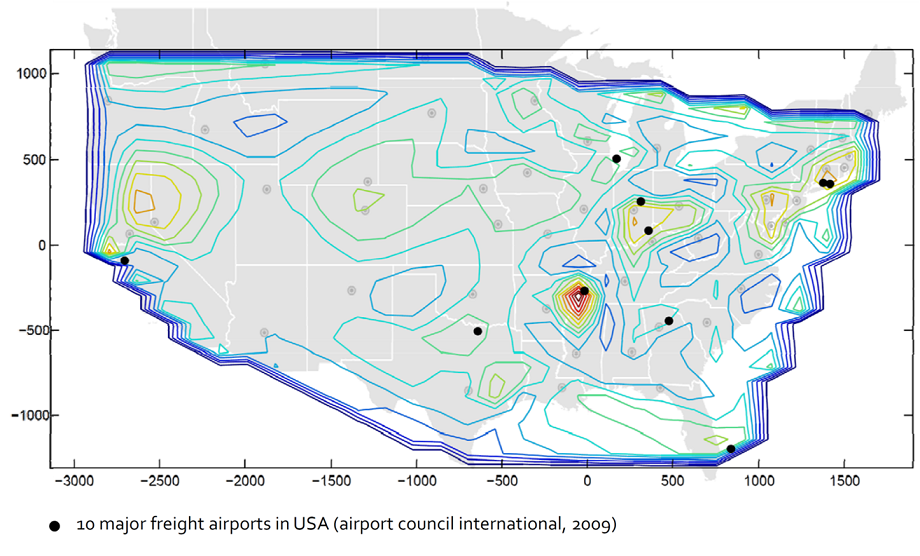}
\caption{density of probability $\mu$ among the US domain, level curves}
\label{Us_res}
\end{figure}

The Figure \ref{Us_res} shows the level curves of the $\mu$ function as result of the real case study in which the coefficient $\eps$ is set $\eps\simeq 10^{-1}$.
 
The black dots of Figure \ref{Us_res} indicate the positions of the 10 major cargo airports in US. The global maxima of the $\mu$ density results very close to the Memphis airport which is the busiest center of airfreight transport and the hub of the FedEx: its airfreight volumes are doubled respect to the other airports. Also the other local maxima are located near the effective position of the other airports: largest errors are appreciated along the boundary areas where also the used projection method presents the largest errors.

\begin{figure} [!]
\centering
\includegraphics[height=8cm]{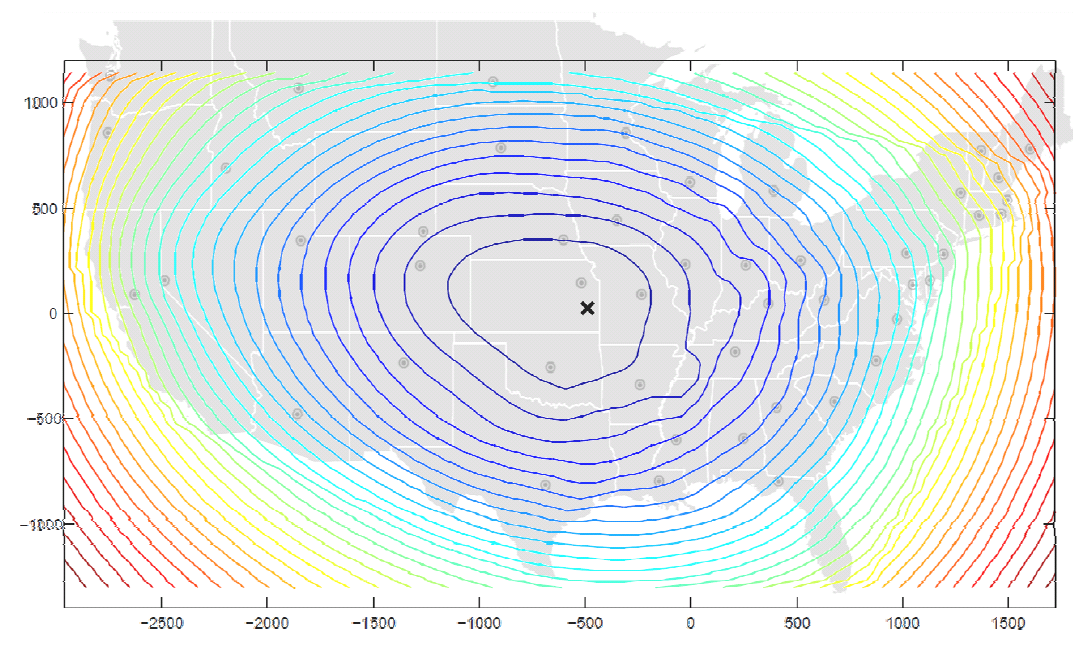}
\caption{Level curves of the functional $H(\delta_{x_0})$}
\label{USfunctional}
\end{figure}

The level curves of the functional $H(\delta_{x_0})$ are plotted in Figure \ref{USfunctional}. Although the exponent $q$ results lower than the unity, the effects of the routing costs tend to predominate the location ones and consequently the importance of initial population $\rho$ is reduced on despite of the distance power relationship.

The position of the minimum point (the ``main hub'') differs from the global maxima of the limit density $\mu$ displayed in Figure \ref{Us_res}. In this case, the difference is due to the value of the coefficient $\eps$ used to determine the limit density $\mu$, which is relatively high so that the results of the functional \eqref{etotnom} and \eqref{ex0} are not coincident.

\section{Conclusions}\label{sconc}

Some considerations can be done. Our initial problem was to determine the optimal position of a certain number $N$ of airports into a domain with location and routing cost condition. This problem is hard if we try to solve it with a direct approach because of his intrinsic complexity. After the modelling phase we concentrate to mass independent routing cost and we caracterized the asymptotic behaviour of the total cost problem. This means that instead of finding the exact position of the N airports, we compute a probability density that represent the ``importance'' of a certain point in the area taken as domain. Moreover, gamma limit result allow us to find the position of the optimal main hub minimizing functional $H(\delta_{x_0})$.

Supported by the examples of the one and two dimensional cases, we observe that these two limit problems are very ``easy'', in terms of computational costs. So instead of looking at the initial problem is more convenient to be reduced at the other two. This makes possible to apply the procedure described to real cases, as in the USA airfreight system.

\bigskip

\ack The work of Giuseppe Buttazzo and Serena Guarino is part of the project 2008K7Z249 {\it``Trasporto ottimo di massa, disuguaglianze geometriche e funzionali e applicazioni''}, financed by the Italian Ministry of Research. The authors wish to thank Giovanni Alberti for some clarifying discussions on the subject.


\bigskip

{\small\parindent0pt
Giuseppe Buttazzo - Dipartimento di Matematica - Universit\`a di Pisa\\
Largo B. Pontecorvo, 5 - 56127 Pisa - ITALY\\
{\tt buttazzo@dm.unipi.it}

\bigskip

Serena Guarino Lo Bianco - Dipartimento di Matematica - Universit\`a di Pisa\\
Largo B. Pontecorvo, 5 - 56127 Pisa - ITALY\\
{\tt sguarino@mail.dm.unipi.it}

\bigskip

Fabrizio Oliviero - Dipartimento di Ingegneria Aerospaziale - Universit\`a di Pisa\\
Via G. Caruso, 8 - 56122 Pisa - ITALY\\
{\tt fabrizio.oliviero@for.unipi.it}}

\end{document}